\documentclass[11pt,a4paper]{article}

\usepackage{microtype}
\usepackage{graphicx}
\usepackage{subfigure}
\usepackage{amsmath,amssymb} %,amsthm}
\usepackage{url,hyperref}
\usepackage{algorithm}
\usepackage[ruled,vlined, linesnumbered, algo2e]{algorithm2e}
\usepackage{algorithmicx,algpseudocode,algcompatible}
\usepackage{mathtools}
\usepackage{multirow}
\usepackage{color}
\usepackage[T1]{fontenc} % Use 8-bit encoding that has 256 glyphs
\usepackage[letterpaper, margin=1.15in]{geometry}
\usepackage{multicol} % Used for the two-column layout of the document
\usepackage[hang, small,labelfont=bf,up,textfont=it,up]{caption} % Custom captions under/above floats in tables or figures
\usepackage{booktabs} % Horizontal rules in tables
\usepackage{float} % Required for tables and figures in the multi-column environment - they need to be placed in specific locations with the [H] (e.g. \begin{table}[H])
\usepackage{hyperref} % For hyperlinks in the PDF

\usepackage{hyperref}

\newcommand{\R}{\mathbb{R}}
\newcommand{\Rext}{\mathbb{R}\cup\{+\infty\}}
\newcommand{\mbf}[1]{\mathbf{#1}}
\newcommand{\mbb}[1]{\mathbb{#1}}
\newcommand{\mcal}[1]{\mathcal{#1}}
\newcommand{\iprods}[1]{\langle{#1}\rangle}
\newcommand{\norms}[1]{\Vert{#1}\Vert}
\newcommand{\sets}[1]{\{{#1}\}}
\newcommand{\BigO}[1]{\mathcal{O}\left({#1}\right)}

\title{Gradient Descent-Type Methods -- Background and Simple Unified Convergence Analysis}

\author{Quoc Tran-Dinh$^{*}$ \textit{and} Marten van Dijk$^{\dagger}$\vspace{0.25ex}\\
\newline {$^{*}$Department of Statistics and Operations Research}\\
\newline {The University of North Carolina at Chapel Hill} \\
\newline {$^{\dagger}$Centrum Wiskunde \& Informatica, Amsterdam, The Netherlands} \\
\newline \textit{Email:} \url{quoctd@email.unc.edu} \textit{and} \url{marten.van.dijk@cwi.nl}}

\date{}
\begin{document}
\maketitle

%%% Abstract.
\begin{abstract}
\normalfont
In this book chapter, we briefly describe the main components that constitute the gradient descent method and its accelerated and stochastic variants.
We aim at explaining these components from a mathematical point of view, including theoretical and practical aspects, but at an elementary level. 
We will focus on basic variants of the gradient descent method and then extend our view to recent variants, especially variance-reduced stochastic gradient schemes (SGD).
Our approach relies on revealing the structures presented inside the problem and the assumptions imposed on the objective function.
Our convergence analysis unifies several known results and relies on a general, but elementary recursive expression.
We have illustrated this analysis on several common schemes.
\end{abstract}

%\begin{keywords}
%\kwd{Gradient descent method}
%\kwd{stochastic gradient method}
%\kwd{convex and nonconvex optimization}
%\kwd{convergence guarantees}
%\kwd{complexity estimates}
%\end{keywords}

%%% 1. Introduction.
\section{Introduction}\label{chapter_XX_intro}
The core problem in many optimization applications such as signal and image processing, engineering, operations research, statistics, and machine learning is the following optimization problem, see, e.g., \cite{Bauschke2011,BenTal2001, Boyd2004,goodfellow2016deep,Hastie2009,Nocedal2006,Sra2020}:
\begin{equation}\label{eq:opt_prob}
\min_{w\in\R^p} F(w),
\end{equation}
where $F : \R^p \to \Rext$ is a given objective or  loss function, and $w$ is a vector of decision variables (also called model parameters).
Depending on the form or structures of the objective function $F$, we obtain different classes of optimization problems.
For instance, the following structures are common in practice.
%%%%%
\begin{itemize}
\item \textbf{Nonsmooth convex optimization.} 
If $F$ is $M$-Lipschitz (i.e. there exists $M > 0$ such that $\vert F(w) - F(w')\vert \leq M\Vert w - w'\Vert$ for all $w, w'\in\R^p$) and convex, but often nonsmooth, then  \eqref{eq:opt_prob} is called a nonsmooth convex minimization.
Note that the $M$-Lipschitz continuity is often imposed for nonsmooth functions such as $F(w) := \norms{w}$ for any norm, or for special smooth functions, e.g., the objective $F(w) := \sum_{i=1}^n\log(1 + \exp(y_iX_i^{\top}w))$ of a logistic regression,  where $(X_i, y_i)$ is given for $i=1,\cdots, n$.
Obviously, the Lipschitz continuity also holds if we consider $F$ to be continuous on a given compact set $\mcal{W}$.

\item \textbf{Smooth and convex optimization.} If $F$ is $L$-smooth (i.e. there exists $L \geq 0$ such that $\Vert \nabla{F}(w) - \nabla{F}(w')\Vert \leq L\Vert w - w'\Vert$ for all $w, w'\in\R^p$) and convex, then  \eqref{eq:opt_prob} is called a smooth and convex minimization.
Examples of $L$-smooth and convex functions are vast.
For example, a least-squares function $F(w) := \frac{1}{2}\norms{X^{\top}w - y}^2$ for a given data matrix $X$ and  an output vector $y$ is $L$-smooth with $L := \norms{XX^{\top}}$. 
The logistic regression function above is also convex and $L$-smooth with $L := \frac{1}{4}\Vert XX^{\top}\Vert$.
However, exponential functions such as $F(w) := \sum_{i=1}^n\exp(X_i^{\top}w)$ or logarithmic functions such as $F(w) := -\sum_{i=1}^n\log(X_i^{\top}w)$ are convex, but not $L$-smooth on their domain, unless we limit their domain on a given compact set, see, e.g., \cite{Tran-Dinh2013a}.

\item \textbf{Smooth and nonconvex optimization.} If $F$ is $L$-smooth and nonconvex, then \eqref{eq:opt_prob} is called a smooth and nonconvex minimization. 
The $L$-smoothness is a key condition required in most gradient-based methods for nonconvex optimization.
Again, this assumption obviously holds if we assume that $F$ is continuously differentiable and then limit the domain of $F$ on a compact set.
But there exists $L$-smooth functions on the entire space $\R^p$.
For instance, $F(w) := \frac{1}{2}w^{\top}Qw + q^{\top}w$ for given symmetric matrix $Q$ and $q\in\R^p$ is $L$-smooth with $L := \norms{Q}$, but not necessarily convex.

\item \textbf{Composite optimization.} If $F(w) := f(w) + g(w)$, where $f$ is usually $L$-smooth and convex/nonconvex, and $g$ is convex and possibly nonsmooth, then \eqref{eq:opt_prob} is called [additive] composite minimization.
This model is ubiquitous in machine learning and statistical learning, where $f$ presents a loss function or a data fidelity term, while $g$ is a regularizer or a penalty term to promote solution structures or to handle constraints. 
Examples can be found, e.g., in \cite{Combettes2011,Pham2019}.
If $g(w)$ is the indicator of a convex set $\mcal{W}$ as $g(w) = 0$ if $w\in\mcal{W}$, and $g(w) = +\infty$, otherwise, then \eqref{eq:opt_prob} covers constrained problem $\min_{w\in\mcal{W}}f(w)$.

\item \textbf{Finite-sum optimization.} If $F(w) := \frac{1}{n}\sum_{i=1}^nF_i(w)$ for some $n \geq 1$, then  \eqref{eq:opt_prob} is called a finite-sum minimization, an empirical risk minimization, or distributed optimization depending on the context.
This structure is presented in most supervised learning tasks, network and distributed optimization, and federated learning.
The most interesting case is when $n\gg 1$.

\item \textbf{Stochastic optimization.} If $F(w) := \mbb{E}[ \mbf{F}(w, \xi)]$, the expectation of a stochastic function $\mbf{F} : \R^p\times\Omega\to\R$, where $(\Omega, \Sigma, \mbb{P})$ is a given probability space, then \eqref{eq:opt_prob} is called a stochastic program \cite{lan2020first,Nemirovskii1983,Shapiro2009}. 
This setting also covers the finite-sum as a special case by setting $\Omega := \sets{1,\cdots, n}$ and $\mbb{P}(\xi=i) = \frac{1}{n}$.
\end{itemize}
Apart from these individual settings, many other combinations between them are also possible; we do not list all of these here.
For example, the combination of composite structure and finite-sum is very common. 

\vspace{1ex}
\noindent\textbf{Existence of solutions.}
We first assume that $F^{\star} := \inf_wF(w)$ is bounded from below, i.e. $F^{\star} > -\infty$ to guarantee the well-definedness of \eqref{eq:opt_prob}.
Many machine learning applications automatically satisfy this condition since the underlying loss function is usually nonnegative. 
One obvious example is the least-squares problem.

Our next question is: \textit{Does \eqref{eq:opt_prob} have an optimal solution?}
To discuss this aspect, we  use a coercive concept from nonlinear analysis \cite{Ekeland1983}.
We say that $F$ is coercive if $\lim_{\norms{w}\to\infty}F(w) = +\infty$.
A common coercive function is $F(w) := \mcal{L}(w) + \frac{\lambda}{2}\norms{w}^2$, where $\mcal{L}$ is $M$-Lipschitz continuous (but not necessarily convex) and $\lambda > 0$.
If $F$ is continuous and coercive, then by the well-known Weierstrass theorem, \eqref{eq:opt_prob} has global optimal solutions $w^{\star}$.
In this case, we denote $F^{\star} := F(w^{\star})$, its optimal value.
If $F$ is nonconvex and differentiable, then we use $w^{\star}$ to denote its stationary points, i.e. $\nabla{F}(w^{\star}) = 0$.
If $F$ is not differentiable, a generalization of stationary points is required \cite{mordukhovich2006variational}.
To keep it simple, we assume throughout this chapter that $F$ is continuously differentiable. 

If $F$ is strongly convex, then it is continuous and coercive and \eqref{eq:opt_prob} has a unique global optimal solution.
For convex problems, our goal is to find an approximate global solution $\hat{w}^{\star}$ of $w^{\star}$ in some sense (see Subsection~\ref{subsec:opt_cert}).
For nonconvex problems, we only expect to compute an approximate stationary point $\hat{w}^{\star}$, which can be a candidate for a local minimizer.
However, we do not attempt to check if it is an approximate local minimizer or not in this chapter.

\vspace{1ex}
\noindent\textbf{Contribution.}
Our contribution can be summarized as follows.
We provide a comprehensive discussion for the main components of the gradient descent method and its variants, including stochastic schemes.
We also propose a unified and simple approach to analyze convergence rates of these algorithms, and demonstrate it through concrete schemes.
This approach can perhaps be extended to analyzing other algorithms, which are not covered in this chapter.
We also discuss some enhanced implementation aspects of the basic algorithms.

\vspace{1ex}
\noindent\textbf{Outline.}
The rest of this chapter is organized as follows.
Section~\ref{chapter_XX_sec_main_body_01} reviews basic components of gradient methods.
Section~\ref{chapter_XX_sec_main_body_02} focuses on stochastic gradient methods, while Section~\ref{chapter_XX_sec_conclusion} makes some concluding remarks and raises a few possible research directions.

\section{Basic components of GD-type methods}\label{chapter_XX_sec_main_body_01}
The gradient descent (GD) method is an iterative method aimed at finding an approximate solution of \eqref{eq:opt_prob}.
This dates back to the works of Cauchy in the 19th century, and has been intensively studied in numerical analysis, including optimization for many decades.
During the last two decades, there has been a great surge in first-order methods, especially gradient-type algorithms, due to applications in signal and image processing, modern statistics, and machine learning. 
In this chapter, we do not attempt to review the literature of GD-type methods, but only focus on summarizing their basic components.

Formally, the gradient descent algorithm  starts from a given initial point $w^0 \in \R^p$ and at each iteration $t \geq 0$, it updates
\begin{equation}\label{eq:gd_main_iter} 
w^{t+1} := \mcal{P}(w^t + \eta_td^t),
\end{equation}
where $w^t$ is the current iterate, $\eta_t > 0$ is called a stepsize or learning rate, $d^t$ is called a search direction, and $\mcal{P}$ is an operator to handle constraints or regularizers;
if this is not needed, one can set $\mcal{P} = \mbb{I}$, the identity operator. 
This method generates a sequence of iterate vectors $\sets{w^t}$ using only first-order information of $F$ (e.g., function values, proximal operators, or [sub]gradients).
Here, we add an operator $\mcal{P}$, which can also be used to handle constraints, regularizers, penalty, or Bregman distance (e.g., mapping between the primal and dual spaces).
Let us discuss each component of the scheme \eqref{eq:gd_main_iter}.

\subsection{Search direction}
The most important component in \eqref{eq:gd_main_iter} is the search direction $d^t$, which determines the type of algorithm such as first-order, second-order, quasi-Newton-type, or stochastic methods.
Let us consider the following possibilities.
\begin{itemize}
\item \textbf{Descent direction.} Assume that $F$ is continuously differentiable, a search direction $d^t$ is called a descent direction at the iterate $w^t$ if $\iprods{\nabla{F}(w^t), d^t} < 0$.
We can even impose a stronger condition, called strictly descent, which is $\iprods{\nabla{F}(w^t), d^t} \leq -c\norms{\nabla{F}(w^t)}^2$ for some $c > 0$.
The name ``descent'' comes from the fact that if we move from $w^t$ along the direction $d^t$ with an appropriate stepsize $\eta_t$, then we have a descent, i.e. $F(w^{t+1}) < F(w^t)$.
If $\mcal{P} = \mbb{I}$, the identity operator, then \eqref{eq:gd_main_iter} reduces to $w^{t+1} := w^t + \eta_td^t$.
By Taylor's expansion of $F$, we have 
\begin{equation*}
F(w^{t+1}) = F(w^t) + \eta_t\iprods{\nabla{F}(w^t), d^t} + o(\eta_t^2\norms{d^t}^2) < F(w^t),
\end{equation*}
for sufficiently small $\eta_t > 0$ due to $\iprods{\nabla{F}(w^t), d^t} < 0$.

\item \textbf{Steepest descent direction.} 
If we take $d^t := -\nabla{F}(w^t)$, then \eqref{eq:gd_main_iter} becomes
\begin{equation}\label{eq:gd_scheme0} 
w^{t+1} := w^t - \eta_t\nabla{F}(w^t),
\end{equation}
and we have $\iprods{\nabla{F}(w^t), d^t} = -\norms{\nabla{F}(w^t)}^2 < 0$ provided that $w^t$ is not a stationary point of $F$.
With this choice of $d^t$, we obtain a gradient descent or also called a steepest descent method.
It actually realizes the most decrease of $F$ at $w^t$ as $\iprods{\nabla{F}(w^t), d^t} \geq -\norms{\nabla{F}(w^t)}$ for any $d^t$ such that $\norms{d^t} = 1$.

\item \textbf{Stochastic gradient direction.}
If we choose $d^t$ to be a stochastic estimator of $\nabla{F}(w^t)$, then we obtain a stochastic approximation (or also called stochastic gradient descent) method.
A stochastic gradient direction is generally not a descent one, i.e., $\iprods{\nabla{F}(w^t), d^t} \not{<} 0$.
Examples of stochastic estimators include standard unbiased estimator $v^t := \nabla{\mbf{F}}(w^t, \xi_t)$ and its mini-batch version $v^t := \frac{1}{\vert\mcal{S}_t\vert}\sum_{\xi\in\mcal{S}_t}\nabla{\mbf{F}}(w^t,\xi)$ for a minibatch $\mcal{S}_t$, and various variance-reduced estimators, see, e.g., \cite{Defazio2014,johnson2013accelerating,nguyen2017sarah,schmidt2017minimizing,Tran-Dinh2019a}.

\item \textbf{Newton and quasi-Newton direction.}
We can go beyond gradient-based methods by incorporating second-order information, or curvature of $F$ as $d^t := -B_t^{-1}\nabla{F}(w^t)$, where $B_t$ is a given symmetric and invertible matrix. 
For instance, if $B_t := \nabla^2{F}(w^t)$, then we obtain a Newton method, while if $B_t$ is an approximation to $\nabla^2{F}(w^t)$, then we obtain a quasi-Newton method.

\item \textbf{Inexact descent direction.}
If we do not evaluate the gradient $\nabla{F}(w^t)$ exactly, but allow some error as $d^t = -(\nabla{F}(w^t) + \delta_t)$ for some Gaussian noise $\delta_t$, then we obtain an inexact or noisy gradient method \cite{ge2015escaping}.
Another example is called sign gradient method, which uses $d^t = -\mathrm{sign}(\nabla{F}(w^t))$, the sign of gradient, see, e.g., \cite{moulay2019properties}.
Inexact Newton-type methods compute $d^t$ by approximately solving $B_td^t = -\nabla{F}(w^t)$ such that $\norms{B_td^t + \nabla{F}(w^t)} \leq c\norms{d^t}$ for $c > 0$.
\end{itemize}

Apart from the above examples, other methods such as [block]-coordinate, incremental gradient, Frank-Wolfe or conditional gradient, proximal-point, prox-linear, Gauss-Newton, extragradient, optimistic gradient, and operator splitting techniques can also be  written in the form \eqref{eq:gd_main_iter} by formulating appropriate search directions $d^t$.
For instance, the proximal point method can be viewed as the gradient method applied to its Moreau's envelope, see, e.g., \cite{Rockafellar2004}.

\subsection{Step-size}
The second important component of \eqref{eq:gd_main_iter} is the step-size $\eta_t$.
In machine learning community, this quantity is called a \textit{learning rate}.
Choosing an appropriate $\eta_t$ is a crucial step that affects the performance of the algorithm.
Classical optimization techniques have proposed several strategies, which are known as globalization strategies, including (i) line-search and its variants, (ii) trust-region, and (iii) filter \cite{Conn2000,Fletcher2002,Nocedal2006}.
Line-search and trust-region strategies have been widely used in numerical optimization, see \cite{Conn2000,Nocedal2006}.
In recent optimization algorithms and machine learning training tasks, we often observe the following techniques.
\begin{itemize}
\item \textbf{Constant learning rate.}
Constant learning rates are usually used to derive convergence rates or complexity bounds due to their simplicity.
The algorithm often performs well with a constant learning rate if the problem is ``easy'', e.g., strongly convex and $L$-smooth, but it becomes poor if the landscape of $F$ is complex such as deep neural networks.
Usually, theoretical analysis gives us a range (i.e. an interval, like $\big(0, \tfrac{2}{L}\big)$ in standard gradient methods) to choose a constant learning rate. 
However, this range could easily be underestimated by using global parameters, and does not capture the desired region of optimal solutions.
In practice, nevertheless, we need to tune this learning rate using different strategies such as grid search or bisection, etc.

\item \textbf{Diminishing learning rate.}
Diminishing learning rates are usually used in subgradient or stochastic gradient methods.
One common diminishing learning rate is $\eta_t := \frac{C}{(t + \beta)^{\nu}}$ for some positive constant $C$, a shifting factor $\beta$, and an order $\nu > 0$.
Depending on the structure assumptions of $F$, we can choose appropriate order $\nu$, e.g., $\nu := \frac{1}{2}$ or $\nu := \frac{1}{3}$.
Other possibility is to choose $\eta_t := \frac{C}{(\lceil t/s\rceil + \beta)^{\nu}}$ for an additional integer $s$ to maintain fixed learning rate in each $s$ iterations.
In stochastic gradient methods, diminishing learning rates are often required if the variance of $d^t$ is nondecreasing (i.e. $d^t$ is not computed from a variance-reduced estimator of $\nabla{F}$).
A diminishing learning rate also determines the convergence rate of the underlying algorithm. 

\item \textbf{Scheduled learning rate.}
In practice, we often use a schedule to tune an appropriate learning rate so that we can achieve the best performance. 
Different ideas have been proposed such as using exponential decay rate, cosine annealing, or even with a varying mini-batch size, see, e.g., \cite{loshchilov10sgdr,tran2021smg}.

\item \textbf{Adaptive learning rate.}
The above strategies of learning rate selection do not take into account the local geometry of the objective function. 
They may perform poorly on ``hard'' problems. 
This motivates the use of adaptive learning rates which exploit local landscape or curvature of the objective function.
The first common strategy is linesearch, which [approximately] solves $\min_{\eta > 0}F(w^t + \eta d^t)$ to find $\eta_t$.
If $F$ is quadratic, then we can compute $\eta_t$ exactly by solving this one-variable minimization problem.
However, most algorithms use inexact line-search such as bisection or golden ratio search.
Another strategy is using a Barzilai-Borwein step-size, e.g., $\eta_t := \norms{w^t - w^{t-1}}/\norms{\nabla{F}(w^t) - \nabla{F}(w^{t-1})}$, which gives an estimation of $\frac{1}{L}$.

Recently, several adaptive methods have been proposed, see, e.g., \cite{Cutkosky2019,Duchi2011,Kingma2014}.
The underlying learning rate is usually given by $\eta_t := C/\sqrt{\sum_{j=0}^t\norms{g_j}^2 + \epsilon}$, where $g_j$ is some gradient estimator at iteration $j$, $C > 0$ is given, and $\epsilon \geq 0$ is a small constant to avoid division by zero and provide numerical stability.

\end{itemize} 
Among the above strategies and tricks for selecting learning rates, one can also compute them in different ways when solving a specific problem, even using a ``trial and error'' method or a combination of the above techniques.
The main goal is to tune a good learning rate for such a problem, but still guarantee the convergence of the algorithm. 

\subsection{Proximal operator}
Many problems covered by \eqref{eq:opt_prob} have constraints or nonsmooth objective terms.
For example, we may have $F(w) = f(w) + g(w)$, where $g$ is nonsmooth.
In this case, we cannot use the full gradient of $F$.
One way to handle the nonsmooth term $g$ is to use proximal operators, and in particular, use the projections if we have simple constraints. 
Mathematically, the proximal operator of a proper and lower semicontinuous function $g$ is defined as
\begin{equation}\label{eq:prox_oper}
\mathrm{prox}_{\gamma g}(w) := \mathrm{arg}\min_{z\in\R^p}\Big\{  \gamma g(z) + \tfrac{1}{2}\norms{ z - w}^2 \Big\}, \quad \gamma > 0.
\end{equation}
Note that under appropriate choices of $\gamma$, the minimization problem in \eqref{eq:prox_oper} is strongly convex, and hence has unique solution, leading to the well-definedness of $\mathrm{prox}_{\gamma g}$.
If $g = \delta_{\mcal{W}}$ as the indicator function of a closed and convex set $\mcal{W}$, i.e. $\delta_{\mcal{W}}(w) = 0$ if $w \in \mcal{W}$, and $\delta_{\mcal{W}}(w) = +\infty$, otherwise, then $\mathrm{prox}_{\gamma g}$ reduces to the projection onto $\mcal{W}$, i.e. $\mathrm{proj}_{\mcal{W}}(w) := \mathrm{arg}\min_{z\in\mcal{W}}\frac{1}{2}\norms{z - w}^2$.
In terms of computation, evaluating $\mathrm{prox}_{\gamma g}$ is generally as hard as solving a [strongly] convex problem.
%However, 
There are at least three ways of evaluating $\mathrm{prox}_{\gamma g}(\cdot)$, which can be sketched as follows.
%--Marten: which ways? Maybe remove this sentence?
\begin{itemize}
\item\textbf{Separable functions.} 
The most obvious case is when $g$ is component-wise separable as $g(w) := \sum_{j=1}^pg_j(w_j)$ (e.g., $g(w) := \norms{w}_1$), then evaluating $\mathrm{prox}_{\gamma g}$ requires solving $p$ one-variable convex minimization problems, which can be done in a closed form. This idea can be extended to block separable functions, e.g., $g(w) := \sum_{i=1}^n\norms{w_{[i]}}_2$, where $\{w_{[i]}\}_{i=1}^n$ are subvectors.
%for $i=1,\cdots, n$.

\item \textbf{Dual approach.} Moreau's identity $\mathrm{prox}_{\gamma g}(w) = w - \gamma \cdot \mathrm{prox}_{g^{\ast}/\gamma}(w/\gamma)$ suggests that we can compute $\mathrm{prox}_{\gamma g}$ from its Fenchel conjugate $g^{\ast}$.
Since many convex functions have simple conjugates such as norms (e.g., $g(w) = \norms{w}_2$) or Lipschitz continuous functions, this approach is more tractable.

\item \textbf{Optimality approach.} If $g$ is differentiable, then we can directly use its optimality condition $\nabla{g}(z) + \gamma^{-1}(z - w) = 0$, and solve it as a nonlinear equation in $z$.
Examples include $-\log\det(X)$ and $\sum_{i=1}^n\log(1 + \exp(y_iX_i^{\top}w))$.
\end{itemize}
Note that the second and third techniques are only used for convex functions, while the first one can be used for nonconvex functions. The number of convex functions $g$ where $\mathrm{prox}_{\gamma g}(\cdot)$ can be computed efficiently is vast, see, e.g., \cite{Bauschke2011,Parikh2013} for more examples and computational techniques. 

\subsection{Momentum}
One way to explain the role of momentum is to use a dynamical system of the form $\ddot{w}(\tau)  + \psi(\tau)\dot{w}(\tau) + \nabla{F}(w(\tau)) = 0$ rooted from Newton's second law, where $\psi(\tau)\dot{w}(\tau)$ presents a friction or a damping factor.
If we discretize this differential equation using $\ddot{w}(\tau) \approx (w^{t+1} - 2w^t + w^{t-1})/h_t^2$ and $\dot{w}(\tau) \approx (w^t - w^{t-1})/h_t$, then we obtain $(w_{t+1} - 2w_t + w_{t-1})/h_t^2 + \psi_t (w_t - w_{t-1})/h_t + \nabla{F}(w_t) = 0$, leading to $w^{t+1} := w^t - h_t^2\nabla{F}(w^t) + (1 - h_t\psi_t)(w^t - w^{t-1})$, see, e.g., \cite{Su2014}.
Therefore, we can specify momentum variants of \eqref{eq:gd_main_iter} when $\mcal{P} = \mbb{I}$ (the identity operator) as follows.
%One can extends the basic scheme \eqref{eq:gd_main_iter} to different directions.
%For example, assume that $\mcal{P} = \mbb{I}$, one can add an inertial term as 
\begin{equation}\label{eq:momentum_scheme}
w^{t+1} := w^t + \eta_td^t + \beta_t(w^t - w^{t-1}),
\end{equation}
where $\beta_t > 0$ is a momentum stepsize. 
The search direction $d^t$ can be evaluated at $w^t$ leading to a so-called heavy ball method \cite{Polyak1964}.
Alternatively, if $d^t$ is evaluated at an intermediate point, e.g., $z^t := w^t + \beta_t(w^t - w^{t-1})$, then we obtain Nesterov's accelerated scheme in the convex case \cite{Nesterov1983}.
This scheme can be written into two steps as 
\begin{equation}\label{eq:nest_scheme}
z^t := w^t + \beta_t(w^t - w^{t-1}), \quad\text{and} \quad w^{k+1} := z^t + \eta_td(z^t),
\end{equation}
where $d(z^t)$ presents the direction $d^t$ evaluated at $z^t$ instead of $w^t$.
Note that momentum terms do not significantly add 
%incur 
computational costs on top of \eqref{eq:gd_main_iter}.
%However, 
Yet, it can accelerate the algorithm in convex cases \cite{Nesterov1983,Nesterov2005c} (see also Subsection~\ref{subsec:conv_analysis}), and possibly in some nonconvex settings, see, e.g., \cite{lee2021fast,tran2022connection}.

\subsection{Dual averaging variant}
The scheme \eqref{eq:gd_main_iter} can be viewed as a forward update, but in convex optimization, dual averaging schemes are also closely related to  \eqref{eq:gd_main_iter}.
Unlike  \eqref{eq:gd_main_iter}, a dual averaging scheme works as follows.
Starting from $w^0\in\R^p$, for $t\geq 0$, we update
\begin{equation}\label{eq:dual_averg}
w^{t+1} := \mathrm{arg}\min_{w}\Big\{ \sum_{j=0}^t\gamma_j\iprods{g^j, w} + \tfrac{1}{2\eta_t}\norms{w - w^0}^2 \Big\}, 
\end{equation}
where $g^j$ are given dual directions (e.g., $g^j := \nabla{F}(w^j)$), $\gamma_j$ are the weights of $g^j$,
%gradients, 
and $\eta_t$ is a given dual stepsize.
In general settings, we can replace $\frac{1}{2}\norms{w - w^0}^2$ by a general Bregman distance $\mcal{D}(w, w^0)$.
If the norm is the Euclidean norm, then we have $w^{t+1} := w^0 - \eta_t\sum_{j=0}^t\gamma_j g^j$.
If $\eta_t = \eta > 0$ is fixed and we choose $g^j := \nabla{F}(w^j)$, then we have $w^{t+1} = w^0 - \eta\sum_{j=0}^t\gamma_jg^j = w^0 - \eta\sum_{j=0}^{t-1}\gamma_jg^j - \eta\gamma_tg^t = w^{t} - \eta\gamma_tg^t$, which is exactly covered by \eqref{eq:gd_main_iter}.
Therefore, for the Euclidean norm $\frac{1}{2}\norms{w - w^0}^2$, the dual averaging scheme \eqref{eq:dual_averg} is identical to the gradient descent scheme $w^{t+1} = w^t - \eta\gamma_tg^t$.
However, under a non-Euclidean norm or a Bregman distance, these methods are different from each other.

\subsection{Structure assumptions}
One main theoretical task when designing a gradient-based algorithm is to establish its convergence.
From a computational perspective, estimating the convergence rate as well as complexity is also critically important.
However, to establish these, we require $F$ to satisfy a set of assumptions.
The following structures are commonly used in optimization modeling and algorithms.
\begin{itemize}
\item \textbf{Lipschitz continuity.} $F$ in \eqref{eq:opt_prob} is said to be $M$-Lipschitz continuous if
\begin{equation}\label{eq:Lipschitz_continuity}
\vert F(w) - F(w')\vert \leq M\Vert w - w'\Vert, \quad \forall w, w'\in\R^p.
\end{equation}
Examples of Lipschitz continuous functions include norms, smoothed approximation of norms (e.g., $F(w) := \sum_{i=1}^p(w_j^2 + \epsilon^2)^{1/2}$ for a small $\epsilon$), or continuous functions with bounded domain.
Note that when $F$ is convex, then $M$-Lipschitz continuity is equivalent to $M$-bounded [sub]gradient, i.e., $\norms{\nabla{F}(w)} \leq M$ for all $w \in \R^p$.
This assumption is usually used in subgradient-type or stochastic gradient-type methods.

\item \textbf{$L$-smoothness.} $F$ is called $L$-smooth if the gradient $\nabla{F}$ of $F$ satisfies 
\begin{equation}\label{eq:L_smoothness}
\Vert \nabla{F}(w) - \nabla{F}(w')\Vert \leq L\Vert w - w'\Vert, \quad \forall w, w'\in\R^p.
\end{equation}
If $w, w'\in\mcal{W}$, for a compact domain $\mcal{W}$ and $F$ is continuously differentiable, then $F$ is $L$-smooth on $\mcal{W}$.
This concept can be extended to an $L$-average smoothness in the finite-sum or stochastic settings.
For instance, if $F(w) := \frac{1}{n}\sum_{i=1}^nF_i(w)$, then we can modify \eqref{eq:L_smoothness} as $\frac{1}{n}\sum_{i=1}^n\norms{\nabla{F}_i(w) - \nabla{F}_i(w')}^2 \leq L^2\norms{w - w'}^2$ for all $w, w'\in\R^p$.
Alternatively, if $F(w) := \mbb{E}[\mbf{F}(w, \xi)]$, then we can use $\mbb{E}[\norms{\nabla\mbf{F}(w) - \nabla\mbf{F}(w')}^2 \mid w, w'] \leq L^2\norms{w - w'}^2$ for all $w, w'\in\R^p$.
These assumptions are usually used in variance reduction SGD methods, see, e.g., \cite{Pham2019,Tran-Dinh2019a}.
Note that other extensions are possible, see, e.g., \cite{lan2020first}.
Verifying the $L$-smoothness is generally not straightforward.
However, if $F(w) := \frac{1}{n}\sum_{i=1}^n\ell_i(X_i^{\top}w - y_i)$ as, e.g., in a generalized linear model, then we can verify the $L$-smoothness of $F$ by verifying the $L$-smoothness of each one-variable function $\ell_i$.
This model is ubiquitous in machine learning. 

One key property of \eqref{eq:L_smoothness} is the following bound:
\begin{equation}\label{eq:L_smoothness2}
\vert F(w') - F(w) - \iprods{\nabla{F}(w), w' - w} \vert \leq \frac{L}{2}\norms{w' - w}^2,
\end{equation}
which shows that $F$ can be globally upper bounded by a convex quadratic function and globally lower bounded by a concave quadratic function.
If, additionally, $F$ is convex, then stronger bounds as well as the co-coerciveness of $\nabla{F}$ can be obtained, see, e.g., \cite{Nesterov2004}.
One can also extend the $L$-smoothness of $F$ to a H\"{o}lder smoothness as $\norms{ \nabla{F}(w) - \nabla{F}(w')}  \leq L\Vert w - w'\Vert^{\nu}$ for some $0 \leq \nu \leq 1$.
This concept unifies both the $L$-smoothness ($\nu = 1$) and the bounded gradient ($\nu = 0$) conditions in one.
It has been used in universal first-order methods for both deterministic and stochastic first-order methods, e.g., \cite{Nesterov2014}.

\item \textbf{Convexity.} $F$ is said to be $\mu$-[strongly] convex if 
\begin{equation}\label{eq:convexity}
F(\hat{w}) \geq F(w) + \iprods{\nabla{F}(w), \hat{w} - w} + \frac{\mu}{2}\norms{\hat{w} - w}^2, \quad \forall w, \hat{w}\in\R^p.
\end{equation}
Here, $\nabla{F}(w)$ can be a gradient or a subgradient of $F$ at $w$.
This inequality shows that $F$ can be lower bounded by either a linear ($\mu=0$) or a quadratic approximation ($\mu \neq 0$).
If $\mu = 0$, then $F$ is just convex or merely convex.
If $\mu > 0$, then $F$ is strongly convex, and $\mu$ is called the strong convexity parameter.
If $\mu < 0$, then $F$ is called weakly convex.
Convexity and strong convexity are key concepts in convex analysis, optimization, and related fields, see, e.g., \cite{Boyd2004,Rockafellar1970}, and we do not further discuss them here.
\end{itemize}
These are three key and also basic assumptions to analyze convergence of \eqref{eq:gd_main_iter} and its variants. 
Nevertheless, other assumptions can also be exploited.
For example, the following conditions are commonly used in different methods.
\begin{itemize}
\item \textbf{Gradient dominance and PL condition.} 
$F$ is called $\sigma$-gradient dominant if $F(w) - F(w^{\star}) \leq \sigma\norms{\nabla{F}(w)}^2$ for all $w\in\R^p$ and $w^{\star}$ is a minimizer of $F$.
Clearly, if $F$ is strongly convex, then it is gradient dominant.
However, there exists nonconvex functions that are gradient dominant.
Note that one can consider local gradient dominance by limiting $w$ in a neighborhood of $w^{\star}$.
We can also extend this concept to different variants.
The gradient dominant condition allows us to obtain a convergence guarantee on the objective residual $F(\hat{w}) - F(w^{\star})$ even in the nonconvex setting. 
Note that this condition is also called Polyak--{\L}ojasiewicz (PL) condition.
These conditions can be used to establish linear convergence or linear-like convergence rates (i.e. linearly converge to a small neighborhood of an optimal solution) \cite{Pham2019,polyak_condition}.

\item \textbf{Uniform convexity and star-convexity.} $F$ is said to be $\mu$-H\"{o}lder uniformly convex of order $\nu \geq 1$ if $F(w') \geq F(w) + \iprods{\nabla{F}(w), w' - w} + \frac{\mu}{\nu}\norms{w' - w}^{\nu}$ for all $w, w' \in \R^p$, see, e.g.,  \cite{zalinescu1983uniformly}.
Clearly, if $\nu = 2$, then we obtain the strong convexity.
If $\nu = 2$ and $w = w^{\star}$, a minimizer of $F$, then $F$ is said to be $\mu$-star strongly convex.
These conditions are often used in gradient-type methods to establish linear convergence rates  \cite{necoara2015linear}.

\item \textbf{Sharpness, quadratic growth, and error bound conditions.}
Assume that  there exist $\gamma > 0$ and $\nu \geq 1$ such that $F(w) - F(w^{\star}) \geq \frac{\gamma}{\nu}\norms{w - w^{\star}}^{\nu}$ for all $w \in \R^p$ and a minimizer $w^{\star}$ of $F$.
If $\nu = 1$, then we say that $F$ is sharped at $w^{\star}$.
If $\nu = 2$, then we say that $F$ has a quadratic growth property.
Clearly, if $F$ is strongly convex, then it has a quadratic growth property.
However, nonconvex functions may still have a quadratic growth property.
This property can be extended to an $\omega$-convexity as in \cite{van2019characterization}.
Another related concept is error bound \cite{luo1993error}, which is defined as $\gamma\norms{\nabla{F}(w)} \geq \norms{w - w^{\star}}$ for some $\gamma > 0$ and all $w\in\R^p$.
Both quadratic growth and error bound conditions can be used to establish [local] linear convergence of gradient-type methods, see, e.g., \cite{drusvyatskiy2018error}.
\end{itemize}
Other properties can be used to analyze convergence of gradient methods such as essential strong convexity, weak strong convexity, restricted secant inequality \cite{necoara2015linear,polyak_condition}, Kurdyka-{\L}ojasiewicz (KL) condition  \cite{bolte2007lojasiewicz}, and Aubin's property \cite{Rockafellar2004}.

\subsection{Optimality certification}\label{subsec:opt_cert}
Finding an exact solution of \eqref{eq:opt_prob} is impractical.
Our goal is to approximate a solution of this problem in some sense.
Let us discuss what we can approximate for \eqref{eq:opt_prob} in both convex and nonconvex problems.

Assume that $w^{\star}$ is a global optimal solution of \eqref{eq:opt_prob} with the optimal value $F^{\star} = F(w^{\star})$ and $\hat{w}$ is an approximate solution produced by an algorithm.
One obvious condition to certify the optimality is to compute the objective residual $F(\hat{w}) - F(w^{\star})$.
We often expect to find $\hat{w}$ such that $F(\hat{w}) - F(w^{\star}) \leq \epsilon$ for a given tolerance $\epsilon > 0$.
This condition is usually used for convex optimization or special classes of nonconvex problems, e.g., under a gradient dominance condition.
The construction of $\hat{w}$ usually relies on two possible ways.
The first one is to simply take the last iterate $w_T$ as $\hat{w} := w^T$, where $w^T$ is the final iterate of the algorithm.
The second option is to form an averaging or a weighted averaging vector as
\begin{equation}\label{eq:averaging_sq}
\hat{w} := \tfrac{1}{T+1}\sum_{t = 0}^Tw^t, \quad \text{or} \quad \hat{w} := \tfrac{1}{S_T}\sum_{t=0}^T\gamma_tw^t, 
\end{equation}
where $\gamma_t > 0$ are given weights (usually related to the stepsize $\eta_t$, but could be different), and $S_T := \sum_{t=0}^T\gamma_t$.
In general, averaging vectors have better theoretical convergence rate guarantees, but they may break desired properties of solutions such as sparsity or low-rankness, etc., compared to the last-iterate $w^T$.
In convex optimization, we often use Jensen's inequality to obtain $F(\hat{w}) - F(w^{\star}) \leq \frac{1}{S_T}\sum_{t=0}^T\gamma_t[F(w^t) - F(w^{\star})]$ for our convergence rate bounds since we obtain a convergence rate bound for the right-hand side.

The second criterion is to use the norm of gradient of $F$, e.g., $\norms{\nabla{F}(\hat{w})}$ or its squared norm.
Note that $\nabla{F}(w^{\star}) = 0$ only provides us stationary points, which are candidates for local minimizers in  nonconvex settings.
Hence, any vector $\hat{w}$ such that $\norms{\nabla{F}(\hat{w})} \leq \epsilon$ for a given tolerance $\epsilon > 0$ only provides us an approximate stationary point $\hat{w}$ of \eqref{eq:opt_prob}. 
To guarantee an approximate local solution, we may add a second-order condition such as $\lambda_{\min}(\nabla^2{F}(\hat{w})) \geq -\hat{\epsilon}$ for some $\hat{\epsilon} > 0$, where $\lambda_{\min}(\nabla^2{F}(\hat{w}))$ is the smallest eigenvalue of $\nabla^2{F}(\hat{w})$.
The construction of $\hat{w}$ in the nonconvex case often relies on the best iterate from $\sets{w^0, \cdots, w^T}$, in the sense that $\norms{\nabla{F}(\hat{w})} = \min_{0\leq t \leq T}\norms{\nabla{F}(w^t)}$.
For nonsmooth optimization, where $F := f + g$, we can use the norm $\norms{G_{\beta}(\hat{w})}$ of gradient mapping $G_{\beta}(\hat{w}) := \beta^{-1}\big(\hat{w} - \mathrm{prox}_{\beta g}(\hat{w} - \beta\nabla{f}(\hat{w}))\big)$ for some $\beta > 0$.
For stochastic optimization, one needs to characterize the optimality condition using expectation $\mbb{E}[\norms{\nabla{F}(\hat{w})}^2]$, $\mbb{E}[F(\hat{w}) - F^{\star}]$, or high probability $\mbb{P}[\norms{\nabla{F}(\hat{w})} \leq \epsilon] \geq 1 - \delta$ or $\mbb{P}[F(\hat{w}) - F^{\star} \leq \epsilon] \geq 1 - \delta$ for a small $\delta\in (0,1)$.

\subsection{Unified convergence analysis}\label{subsec:conv_analysis}
Let us first present our general and unified convergence analysis approach and then illustrate it through three different methods.

\vspace{1ex}
\noindent\textbf{(a) General approach.}
Most convergence analysis of first-order methods of the form \eqref{eq:gd_main_iter} relies on the following recursive inequality often generated by two or three consecutive iterates:
\begin{equation}\label{eq:recursive_est}
D_{t+1} + \Delta_t \leq \omega_t \cdot D_t + E_t,
\end{equation}
where $D_t$, $\Delta_t$, and $E_t$ are nonnegative quantities, and $\omega_t \in (0, 1]$ is a contraction factor.
Very often these quantities depend on two consecutive iterates $w^t$ and $w^{t+1}$, but sometimes they also depend on $w^{t-1}$.
The error term $E_t$ usually satisfies $\sum_{t=0}^{\infty}E_t < +\infty$.
Moreover, we often have $\omega_t = 1$ or $\omega_t \to 1$ for sublinear rates, and a fixed $\omega_t = \omega \in (0, 1)$ for linear rates.
The quantity $D_t$ can be referred to as a potential or Lyapunov function. 
There is no general and universal method to construct $D_t$, but for gradient-type methods, it is usually either $\norms{w^t - w^{\star}}^2$, $\norms{w^t - w^{t-1}}^2$, $F(w^t) - F^{\star}$, $\norms{\nabla{F}(w^t)}^2$ (in Euclidean or weighted norms), or a combination of these terms.
Clearly, if $E_t = 0$, then $\sets{D_t}$ is nonincreasing, showing a descent property of $D_t$.
However, if $E_t > 0$, then we no longer have a descent property of $D_t$, which is usually the case in SGD or subgradient methods.
There are two cases.
\begin{itemize}
\item[]\textbf{Case 1.} If $D_t$ contains an optimality measure, e.g., $S_t[F(w^t) - F^{\star}]$, then we can show that $F(w^t) - F^{\star} \leq \frac{C}{S_t}$ for the last iterate $w^t$, where  $C$ is a constant depending on $w^0$ and possibly on  $w^{\star}$ or $F^{\star}$.
\item[]\textbf{Case 2.} If $\Delta_t$ contains an optimality measure, e.g., $\gamma_t\norms{\nabla{F}(w^t)}^2$, then we can show that $\frac{1}{S_T}\sum_{t=0}^T\gamma_t\norms{\nabla{F}(w^t)}^2 \leq \frac{C}{S_T}$ for some constant $C$ and $S_T := \sum_{t=0}^T\gamma_t$.
\end{itemize}
The recursive estimate \eqref{eq:recursive_est} can be used to prove the convergence of different gradient-type methods, including standard and accelerated algorithms.
Let us illustrate how to obtain \eqref{eq:recursive_est} for some common schemes.

\vspace{1ex}
\noindent\textbf{(b)~Subgradient method.}
Let us consider the classical [sub]gradient method to minimize $F(w)$ as $w^{t+1} = w^t - \eta_t\nabla{F}(w^t)$, which is a special case of \eqref{eq:gd_main_iter}, where $\nabla{F}(w^t)$ is a [sub]gradient of $F$ at $w^t$.
Then, for any $w \in \R^p$, we have 
\begin{equation}\label{eq:proof1a}
\arraycolsep=0.2em
\begin{array}{lcl}
\eta_t\iprods{\nabla{F}(w^t), w^t - w}  =  \tfrac{1}{2}\norms{w^t - w}^2 - \tfrac{1}{2}\norms{w^{t+1} - w}^2 + \tfrac{\eta_t^2}{2}\norms{\nabla{F}(w^t)}^2.
\end{array}
\end{equation}
If $F$ is convex, then $\iprods{\nabla{F}(w^t), w^t - w} \geq F(w^t) - F(w)$.
Combining this inequality and \eqref{eq:proof1a}, we obtain
\begin{equation}\label{eq:proof1b}
\arraycolsep=0.2em
\begin{array}{lcl}
\underbrace{\tfrac{1}{2}\norms{w^{t+1} - w}^2}_{D_{t+1}} + \underbrace{\eta_t[F(w^t) - F(w)] }_{\Delta_t} \leq  \underbrace{\tfrac{1}{2}\norms{w^t - w}^2}_{D_t}  + \underbrace{\tfrac{\eta_t^2}{2}\norms{\nabla{F}(w^t)}^2}_{E_t}.
\end{array}
\end{equation}
This inequality is exactly in the form \eqref{eq:recursive_est} with $\omega_t = 1$.
To guarantee convergence, we need to take $w = w^{\star}$ as a solution of \eqref{eq:opt_prob} and assume that $\norms{\nabla{F}(w^t)} \leq M$.
Then, \eqref{eq:proof1b} implies that $\Delta_t \leq D_t - D_{t+1} + E_t$.
By induction, we have $\sum_{t=0}^T\Delta_t \leq D_0 - D_{T+1} + \sum_{t=0}^TE_t \leq D_0 + \sum_{t=0}^TE_t$.
Therefore, we obtain 
\begin{equation*}
F(\hat{w}) - F(w^{\star}) \leq \frac{1}{S_T}\sum_{t=0}^T\eta_t [F(w^t) - F(w^{\star})] \leq \frac{1}{2S_T}\norms{w^0 - w^{\star}}^2 + \frac{M^2}{2S_T}\sum_{t=0}^T\eta_t^2, 
\end{equation*}
where $S_T := \sum_{t=0}^T\eta_t$ and $\hat{w} := \frac{1}{S_T}\sum_{t=0}^T\eta_tw^t$ as computed by \eqref{eq:averaging_sq}.
To obtain a convergence rate bound, we require  $\sum_{t=0}^{\infty}\eta_t^2 <+\infty$ and $S_T \to S_{\infty} = \sum_{t=0}^{\infty}\eta_t = \infty$.
These are exactly the conditions to guarantee the convergence of [sub]gradient methods, see, e.g., \cite{Boyd2003}.

\vspace{1ex}
\noindent\textbf{(c) Gradient descent method for nonconvex problems.}
If we assume that $F$ is only $L$-smooth and not necessarily convex, then using \eqref{eq:L_smoothness2} with $w := w^t$ and $w' := w^{t+1}= w^t - \eta_t\nabla{F}(w^t)$ we have
\begin{equation}\label{eq:proof1c}
\arraycolsep=0.2em
\begin{array}{lcl}
F(w^{t+1}) & \leq & F(w^t) + \iprods{\nabla{F}(w^t), w^{t+1} - w^t} + \tfrac{L}{2}\norms{w^{t+1} - w^t}^2 \vspace{1ex}\\
%&= & F(w^t) + \iprods{\nabla{F}(w^t), w - w^t} + \vspace{1ex}\\
&= & F(w^t) -\eta_t\big(1 - \frac{L\eta_t}{2}\big)\norms{\nabla{F}(w^t)}^2.
\end{array}
\end{equation}
By adding $-F^{\star}$, where $F^{\star} := \inf_wF(w) > -\infty$ (our assumption), to both sides and rearranging the result, the inequality \eqref{eq:proof1c} leads to 
\begin{equation*}
\arraycolsep=0.2em
\begin{array}{lcl}
\underbrace{F(w^{t+1}) - F^{\star}}_{D_{t+1}} + \underbrace{\eta_t\big(1 - \tfrac{L\eta_t}{2}\big)\norms{\nabla{F}(w^t)}^2}_{\Delta_t} & \leq & \underbrace{F(w^t) - F^{\star}}_{D_t}.
\end{array}
\end{equation*}
This is exactly in the form \eqref{eq:recursive_est} with $\omega_t = 1$ and $E_t = 0$.
Without any further assumption, we have $\Delta_t \leq D_t - D_{t+1}$, and by induction, we get $\sum_{t=0}^T\Delta_t \leq D_0 - D_{T+1} \leq D_0$, leading to
\begin{equation*}
\arraycolsep=0.2em
%\begin{array}{lcl}
\min_{0 \leq t\leq T}\norms{\nabla{F}(w^t)}^2 \leq \frac{1}{S_T}\sum_{t=0}^T\gamma_t\norms{\nabla{F}(w^t)}^2 \leq \frac{F(w^0) - F^{\star}}{S_T},
%\end{array}
\end{equation*}
where $\gamma_t = \eta_t(1 - \tfrac{L\eta_t}2{})$ and  $S_T = \sum_{t=0}^T\gamma_t$, provided that $0 < \eta_t < \frac{2}{L}$.
This result allows us to certify the best-iterate convergence rate of the algorithm to a stationary point of \eqref{eq:opt_prob}.

\vspace{1ex}
\noindent\textbf{(d)~Gradient descent method for smooth and convex problems.}
Assume that $F$ is convex and $L$-smooth.
Let us choose $\eta_t := \frac{1}{L}$ in \eqref{eq:gd_main_iter} to get $w^{t+1} := w^t - \frac{1}{L}\nabla{F}(w^t)$.
Then, from \eqref{eq:proof1a} and \eqref{eq:proof1c}, and the convexity of $F$, we have
\begin{equation*} 
\arraycolsep=0.2em
\left\{\begin{array}{lcl}
\tfrac{L}{2}\norms{w^{t+1} - w^{\star}}^2 + \iprods{\nabla{F}(w^t), w^t - w^{\star}}  &= &  \tfrac{L}{2}\norms{w^t - w^{\star}}^2 + \tfrac{1}{2L}\norms{\nabla{F}(w^t)}^2, \vspace{1ex}\\
(t+1)[ F(w^{t+1}) - F(w^{\star})] +   \frac{(t+1)}{2L}\norms{\nabla{F}(w^t)}^2 &\leq & (t+1)[F(w^t) - F(w^{\star})], \vspace{1ex}\\
F(w^t) - F(w^{\star}) & \leq & \iprods{\nabla{F}(w^t), w^t - w^{\star}}.
\end{array}\right.
\end{equation*}
By summing up these three inequalities and canceling terms, we obtain 
\begin{equation*} 
\arraycolsep=0.2em
\begin{array}{lcl}
\tfrac{L}{2}\norms{w^{t+1} - w^{\star}}^2 + (t+1)[ F(w^{t+1}) - F(w^{\star})] & + & \frac{t}{2L}\norms{\nabla{F}(w^t)}^2   \leq   \tfrac{L}{2}\norms{w^t - w^{\star}}^2 \vspace{1ex}\\
&& + {~} t [F(w^t) - F(w^{\star})],
\end{array}
\end{equation*}
This is exactly \eqref{eq:recursive_est} with $D_t :=  \tfrac{L}{2}\norms{w^t - w^{\star}}^2 + t [F(w^t) - F(w^{\star})]$, $\Delta_t :=  \frac{t}{2L}\norms{\nabla{F}(w^t)}^2$, $E_t = 0$, and $\omega_t = 1$.
This recursive estimate implies $D_{t+1}\leq D_0$, and therefore, using the definition of $D_{t+1}$ and dropping $\tfrac{L}{2}\norms{w^{t+1} - w^{\star}}^2$, we get % \textcolor{red}{replaced $t+1$ by $t$}
\begin{equation*}
F(w^{t+1}) - F(w^{\star}) \leq \frac{D_0}{t+1} = \frac{L\norms{w^0 - w^{\star}}^2}{2(t+1)}, 
\end{equation*}
which shows a $\BigO{1/t}$-last-iterate convergence rate on $w^t$.
It also implies that $\sum_{t=0}^T t\norms{\nabla{F}(w^t)}^2 \leq L^2\norms{w^0 - w^{\star}}^2$ (by using $\sum_{t=0}^T\Delta_t \leq  D_0$) 
%\textcolor{red}{replaced $L$ by $L^2$} 
and $\norms{w^t - w^{\star}} \leq \norms{w^0 - w^{\star}}$ (by using $D_t\leq D_0$) for all $t\geq 0$.

\vspace{1ex}
\noindent\textbf{(e)~Accelerated gradient method for smooth and convex problems.}
Our last illustration 
%is the following 
follows Nesterov's accelerated gradient scheme:
\begin{equation}\label{eq:acc_gd_scheme}
z^t := w^t + \beta_t(w^t - w^{t-1}) \quad \text{and} \quad w^{t+1} := z^t - \tfrac{1}{L}\nabla{F}(z^t),
\end{equation}
where $\beta_t = \frac{\theta_{t-1} - 1}{\theta_{t}}$ for $\theta_t \geq 1$ such that $\theta_t(\theta_t - 1) \leq \theta_{t-1}^2$ with $\theta_0 := 1$.
This is an accelerated variant of \eqref{eq:gd_main_iter} with the momentum $\beta_t(w^t - w^{t-1})$.
It is well-known \cite{Nesterov1983} that, after a few elementary transformations, \eqref{eq:acc_gd_scheme} can be written as 
\begin{equation*}
z^t := (1-\tfrac{1}{\theta_t})w^t + \tfrac{1}{\theta_t}u^t,  \ w^{t+1} := z^t - \tfrac{1}{L}\nabla{F}(z^t), \ \text{and} \ u^{t+1} = u^t - \tfrac{\theta_t}{L}\nabla{F}(z^t).
\end{equation*}
Let $v^t := (1-\tfrac{1}{\theta_t})w^t + \tfrac{1}{\theta_t}w^{\star}$.
Then, $z^t - v^t  =  \tfrac{1}{\theta_t}(u^t - w^{\star})$.
Moreover, by convexity of $F$, we have  $F(z^t)   \leq F(v^t)  + \iprods{\nabla{F}(z^t), z^t - v^t} \leq  (1-\tfrac{1}{\theta_t})F(w^t) + \tfrac{1}{\theta_t}F(w^{\star})  + \tfrac{1}{\theta_t}\iprods{\nabla{F}(z^t),  u^t - w^{\star}}$.
Hence, multiplying  both sides by $\theta_t^2$, we obtain
\begin{equation*} 
\arraycolsep=0.2em
\begin{array}{lcl}
\theta_t^2 [F(z^t) - F(w^{\star})] &\leq & \theta_t(\theta_t - 1)[F(w^t) - F(w^{\star})] +  \theta_t\iprods{\nabla{F}(z^t),  u^t - w^{\star}}.
\end{array}
\end{equation*}
Similar to the proof of \eqref{eq:proof1a} and \eqref{eq:proof1c}, respectively we have
\begin{equation*} 
\arraycolsep=0.2em
\left\{\begin{array}{ll}
& \tfrac{L}{2}\norms{u^{t+1} - w^{\star}}^2 + \theta_t\iprods{\nabla{F}(z^t), u^t - w^{\star}}  =   \tfrac{L}{2}\norms{u^t - w^{\star}}^2 + \tfrac{\theta^2_t}{2L}\norms{\nabla{F}(z^t)}^2, \vspace{1ex}\\
& \theta_t^2[F(w^{t+1}) - F(w^{\star})] +   \frac{\theta_t^2}{2L}\norms{\nabla{F}(z^t)}^2 \leq  \theta_t^2[F(z^t) - F(w^{\star})].
\end{array}\right.
\end{equation*}
Summing up the last three inequalities, we obtain 
\begin{equation}\label{eq:proof3} 
\hspace{-0ex}
\arraycolsep=0.1em
%\begin{array}{lcl}
\begin{aligned}
& \underbrace{\theta_t^2[F(w^{t+1}) - F(w^{\star})]  +   \tfrac{L}{2}\norms{u^{t+1} - w^{\star}}^2}_{D_{t+1}} + \underbrace{(\theta_{t-1}^2 - \theta_t(\theta_t - 1)) [F(w^t) - F(w^{\star})]}_{\Delta_t} \vspace{-1ex}\\
&\hspace{4ex} \leq  \underbrace{\theta_{t-1}^2[F(w^t) - F(w^{\star})] + \tfrac{L}{2}\norms{u^t - w^{\star}}^2}_{D_t},
%\end{array}
\vspace{-1ex}
\end{aligned}
\hspace{-10ex}
\end{equation}
which is exactly \eqref{eq:recursive_est} with $E_t = 0$ and $\omega_t = 1$, provided that $\theta_{t-1}^2 - \theta_t(\theta_t - 1) \geq 0$ (note that $\theta_0 = 1$ and $\theta_{-1} = \frac{1}{2}$ satisfy this condition). 
The recursive estimate \eqref{eq:proof3} implies that $D_t \leq D_0$, leading to 
\begin{equation*}
F(w^t) - F(w^{\star}) \leq \frac{D_0}{\theta_{t-1}^2} = \frac{L}{2\theta_{t-1}^2}\norms{w^0 - w^{\star}}^2.
\end{equation*}
In particular, if we choose $\theta_{t-1} := \frac{t+1}{2}$, then $\theta_{t-1}^2 = \frac{(t+1)^2}{4} \geq \theta_t(\theta_t - 1) = \frac{t(t+1)}{4}$, then we get the last-iterate convergence guarantee $F(w^t) - F(w^{\star}) \leq \frac{2L\norms{w^0 -w^{\star}}^2}{(t+1)^2}$.

We have illustrated how to employ the unified recursive expression \eqref{eq:recursive_est} to analyze four different deterministic gradient-type algorithms.
It  provides a simple approach with a few lines  to derive convergence rate analysis compared to classical techniques in the literature.  
We believe that this approach can be extended to analyze other methods that have not been listed here.

\subsection{Convergence rates and complexity analysis}
Classical optimization literature often characterizes asymptotic convergence or linear convergence rates of the underlying algorithm, while sublinear rates or oracle complexity are largely elusive, see, e.g.,  \cite{Bertsekas1999,Fiacco1987,Fletcher1987,kelley1999iterative,Luenberger2007,Polak1971}. 
Sublinear convergence rates have been widely studied in convex optimization methods, see, e.g., \cite{Nesterov2004}, while oracle complexity analysis was formally studied in \cite{Nemirovskii1983}.
Recently, these topics have gained in popularity due to applications to large-scale problems in modern signal and image processing, machine learning, and statistical learning \cite{bubeck2014theory,Hastie2009,Wright2009}.
Let us discuss these concepts in detail here.

\vspace{1ex}
\noindent\textbf{(a) Convergence rates.}
A convergence rate characterizes the progress of the optimality measure (e.g., the objective residual $F(\hat{w}^t) - F^{\star}$, the squared distance to solution $\norms{\hat{w}^t - w^{\star}}^2$, or the squared norm of gradient $\norms{\nabla{F}(\hat{w}^t)}^2$) w.r.t. the iteration $t$, where $\hat{w}^t$ is an approximate solution.
For example, in the gradient method for smooth and convex problems, we have $F(w^{t+1}) - F^{\star} \leq \tfrac{L\norms{w^0 - w^{\star}}^2}{2(t+1)}$ showing that the objective residual $F(w^t) - F^{\star}$ decreases with a speed of at least $\frac{1}{t}$, which we write $F(w^t) - F^{\star} = \BigO{1/t}$.
We can also write $F(w^{t+1}) - F^{\star} = \BigO{\frac{R_0^2L}{t+1}}$ for $R_0 := \norms{w^0 - w^{\star}}$ to show the dependence of the rate on $L$ and $R_0$.

Note that we generally attempt to establish an upper bound rate, but can also show that this upper bound matches the lower bound rate (up to a constant factor) for certain class of algorithms under a given set of assumptions on \eqref{eq:opt_prob}, see, e.g., \cite{Nesterov2004}.
For gradient-type methods, the optimal convergence rates under only convexity and $L$-smoothness is $\BigO{1/t^2}$, which is guaranteed by Nesterov's optimal methods. 
For nonconvex problems, gradient-type methods only achieve a $\BigO{1/t}$ rate on $\norms{\nabla{F}(\hat{w}^t)}^2$ under  $L$-smoothness.
Linear convergence rates can be achieved with additional assumptions such as strong convexity, error bound, quadratic growth, or PL condition.
However, we do not further discuss these variants in this paper.

\vspace{1ex}
\noindent\textbf{(b) Complexity.}
The concept of complexity comes from theoretical computer science, but is widely used in computational mathematics, and in particular, in optimization.
Formal definitions of complexity can be found, e.g., in \cite{Nemirovskii1983,Nesterov2004}.
We distinguish two types of complexity for our gradient-type methods: iteration-complexity (or analytical complexity), and computational complexity (or sometimes called arithmetic complexity, or work complexity) \cite{Nesterov2004}.
In gradient-type methods, the overall computational complexity is generally  dominated by the oracle complexity, which characterizes the total number of function and/or gradient evaluations required for finding an approximate solution. 
%\textcolor{red}{Do we need the next sentence?} 
We notice that, we  overload the concept \textit{oracle}, which is formally defined, e.g., in \cite{Nesterov2004}.
Mathematically, the oracle complexity of $T$ iterations of an algorithm (in our context) is defined as follows:
\begin{equation}\label{eq:oracle_complexity}
\textrm{Oracle complexity} := \sum_{t=0}^T \textrm{Per-iteration complexity at iteration $t$},
\end{equation}
The \textbf{per-iteration complexity} characterizes the workload (e.g., the number of gradient evaluations) at each iteration. 
At each iteration, we often count the most dominated computation steps such as gradient evaluations, function evaluations, proximal operations, projections, matrix-vector multiplications, or Hessian-vector multiplications.
If this per-iteration complexity is fixed, then we have
\begin{equation*}
\textrm{Oracle complexity} = \textrm{Number of iterations} \times \textrm{Per-iteration complexity}.
\end{equation*}
For example, for the standard gradient descent method for smooth and convex problems, the per-iteration complexity is $\BigO{1}$, i.e. requires one gradient evaluation, leading to oracle complexity $\BigO{\frac{1}{\epsilon}}$ in order to obtain $w^t$ such that $F(w^t) - F^{\star} \leq \epsilon$.
Indeed, from the convergence bound $F(w^t) - F^{\star} \leq \frac{L\norms{w^0 - w^{\star}}^2}{2t}$ we infer that $F(w^t) - F^{\star} \leq \epsilon$ is implied by $\frac{L\norms{w^0 - w^{\star}}^2}{2t} \leq \epsilon$, leading to $t \geq \big\lceil \frac{L\norms{w^0 - w^{\star}}^2}{2\epsilon}\big\rceil$.
Hence, we need at most $t_{\max} := \big\lceil \frac{L\norms{w^0 - w^{\star}}^2}{2\epsilon}\big\rceil = \BigO{1/\epsilon}$ iterations, leading to $\BigO{1/\epsilon}$ gradient evaluations.

\subsection{Initial point, warm-start, and restart}
For convex algorithms, which can converge to a global minimizer $w^{\star}$ starting from any initial point $w^0$, the choice of $w^0$ will affect the number of iterations as the term $\norms{w^0 - w^{\star}}^2$ for any solution $w^{\star}$ appears in the bound of the convergence guarantee, e.g., $F(w^T) - F^{\star} \leq \BigO{\frac{L\norms{w^0 - w^{\star}}^2}{T^{\nu}}}$ for $\nu = 1$ or $\nu = 2$. 
Clearly, if $w^0$ is close to $w^{\star}$, then the number of iterations $T$ is small.

For nonconvex algorithms, initialization plays a crucial role since different initial points $w^0$ may make the algorithm converge to different approximate stationary points $w^{\star}$, and their quality 
%of them 
is different.
Stationary points are candidates for local minimizers, but some may give us maximizers or saddle points.
%Although 
If we do get a local minimizer, then it may still be a bad one, 
which is 
far from any global minimizer or which gives us a bad prediction error in machine learning.

A warm-start strategy 
%we mention here 
uses the output of the previous run or the previous iteration to initialize the algorithm at the current stage or iteration.
It is based on the idea that the previous run already gives us a good approximation of the desired solution.
Initializing from this point may hope to quickly converge to the target optimal solution.
Warm-start is widely used in sequential iterative (e.g., sequential quadratic programming) or online learning methods.

A restarting strategy is often 
%mentions
used in the case where the algorithm makes undesired progress and needs to be restarted.
This idea has been used in accelerated gradient methods, where the objective function increases after significant decrease, causing oscillated behaviors \cite{Giselsson2014,Odonoghue2012}. 
Restarting is often combined with a warm-start and an appropriate condition to obtain good performance.
Some theoretical analysis and practical discussion of restarting strategies can be found, e.g., in \cite{fercoq2016restarting,Giselsson2014,Odonoghue2012,Su2014}.

%%%%%%%%%%%%%%%%%%%%%%%%%%%%%%%%%%%%%%%%%%%%%
\section{Stochastic Gradient Descent Methods}\label{chapter_XX_sec_main_body_02}
Let us further extend our discussion from deterministic to stochastic methods for solving \eqref{eq:opt_prob} when $F$ is a finite-sum or an expectation function.
The stochastic approximation (SA) method was initially proposed by Robbins and Monro in 1950s \cite{RM1951}. It has become extremely popular in the last decades as it has been widely used in machine learning and data science, see, e.g., \cite{Bottou1998,Bottou2018,Sra2020}.

\subsection{The algorithmic template}
In this section, we only focus on the standard stochastic optimization and discuss two types of methods: classical SGD and variance-reduced SGD.
More specifically, we focus on $F(w) := \mathbb{E}[\mbf{F}(w, \xi)]$ in \eqref{eq:opt_prob}, which can be written as
\begin{equation}\label{eq:opt_prob2}
\min_{w\in\R^p} \Big\{ F(w) :=  \mathbb{E}[\mbf{F}(w, \xi)] \Big\},
\end{equation}
where $\xi$ is a random vector defined on a given probability space $(\Omega, \Sigma, \mbb{P})$.

Many stochastic gradient-based methods for solving \eqref{eq:opt_prob2} can be described as in Algorithm~\ref{alg:A1}.
\begin{algorithm}[hpt!]\caption{(Unified Stochastic Gradient (SGD) Method)}\label{alg:A1}
\normalsize
\begin{algorithmic}[1]
\State\label{step:i0}{\bfseries Initialization:} Choose an initial point $\hat{w}^0$ in $\R^p$.
\State\hspace{0ex}{\bfseries For $s=0$ to $S-1$, perform:}
\State\hspace{3ex}Evaluate a snapshot estimator $\hat{v}^s$ of $\nabla{F}(\hat{w}^{s})$ and set $w^{s,0} = \hat{w}^s$;
\State\hspace{3ex}{\bfseries For $t = 0$ to $T_s-1$, update:}
\State\hspace{5ex}Sample a subset of examples $\mcal{S}_{s,t}$;
\State\hspace{5ex}Construct an estimator $v^{s,t}$ of $\nabla{F}(w^{s,t})$ using $\mcal{S}_{s,t}$ and $\hat{v}^s$;
\State\hspace{5ex}Update $w^{s,t+1} := \mcal{P}(w^{s,t} - \eta_{s,t}v^{s,t})$;
\State\hspace{3ex}{\bfseries End of Iterations}
\State\hspace{3ex}Form a new snapshot point $\hat{w}^{s+1}$ from $\sets{w^{s,0}, \cdots, w^{s,T_s}}$.
\State\hspace{0ex}{\bfseries End of Stages}
\State\hspace{0ex}{\bfseries Output:} Return $\hat{w}$ from the available iterates.
\end{algorithmic}
\end{algorithm}
Here, Algorithm~\ref{alg:A1} only presents a pure stochastic gradient scheme with a possible variance-reduction step, but without momentum or accelerated steps.
The operator $\mcal{P}$ presents a projection to handle constraints if required, or to add a compression.
However, if it is not specified, then we assume that $\mcal{P}(z) = z$, the identity operator.
Note that Algorithm~\ref{alg:A1} is a double-loop algorithm, where the inner loop carries out SGD updates, while the outer loop performs stage-wise updates, which can be expressed in an epoch-wise fashion or as a restarting mechanism.

If $S=0$, then Algorithm~\ref{alg:A1} reduces to a single-loop method.
If $S >1$, then we can also transform Algorithm~\ref{alg:A1} into a single loop with ``IF'' statement and using the iteration counter $k := \sum_{i=0}^{s-1}T_i + t$.
If $T_s := T$ is fixed, then $k := (s-1)T + t$.
This transformation allows us to inject Bernoulli's rule for the ``IF'' statement instead of deterministic rules.
Such a modification has been implemented in Loopless-SVRG  and Loopless-SARAH schemes, see, e.g., \cite{kovalev2019don,Li2019}. 

\subsection{SGD estimators}\label{subsec:SGD_estimators}
The main component of Algorithm~\ref{alg:A1} is the estimator $v^t$ of $\nabla{F}(w^t)$.
Let us review some important estimators widely used in optimization and related fields.

\vspace{1ex}
\noindent\textbf{(a) Classical SGD and mini-batch estimators.}
Clearly, if $S = 0$, then we can simply drop the superscript $s$ in Algorithm~\ref{alg:A1}, and write the main update as 
\begin{equation}\label{eq:sgd_update}
w^{t+1} := w^t - \eta_tv^t,
\end{equation}
which is in the form \eqref{eq:gd_main_iter} with $d^t = -v^t$ being a stochastic estimator of $\nabla{F}(w^t)$.

In classical SGD, we often generate $v^t$ as an unbiased estimator of $\nabla{F}(w^t)$ with bounded variance, i.e.:
\begin{equation}\label{eq:unbiased_estimator}
\mbb{E}[v^t \mid \mcal{F}_t] = \nabla{F}(w^t) \quad \text{and} \quad \mbb{E}[\norms{v^t}^2 \mid \mcal{F}_t] \leq M^2,
\end{equation}
for given $M \geq 0$, where $\mcal{F}_t$ is the smallest $\sigma$-algebra generated by $\sets{\mcal{S}_0, \cdots, \mcal{S}_t}$ and $\mbb{E}[\cdot \mid\mcal{F}_t]$ is the conditional expectation.

\vspace{1ex}
\noindent\textbf{(b) Variance-reduced SGD estimators.}
There exists a number of variance-reduced methods, which are based on different estimators of $\nabla{F}(w)$.
We only focus on some of them.
For simplicity, we drop the stage superscript ``s''.

The first one is SVRG \cite{johnson2013accelerating}, which generates $v^{t}$ as
\begin{equation}\label{eq:svrg_estimator}
v^{t} := \hat{v}^{s} + [\nabla{\mbf{F}}(w^{t}, \mcal{S}_{t}) - \nabla{\mbf{F}}(\hat{w}^s, \mcal{S}_{t})],
\end{equation}
where $\nabla{\mbf{F}}(w^{t}, \mcal{S}_{t}) := \tfrac{1}{b_t}\sum_{\xi_{t}\in\mcal{S}_{t}}\nabla{\mbf{F}}(w^{t}, \xi_{t})$ and $b_t := \vert \mcal{S}_{t}\vert$.
Then, one can show that
\begin{equation*}
\mbb{E}_{\mcal{S}_{t}}[ v^{t} ] = \nabla{F}(w^{t}) \quad\text{and}\quad \mbb{E}_{\mcal{S}_{t}}[ \norms{ v^{t} - \nabla{F}(w^{t})}^2 ] \leq \sigma_{t}^2, 
\end{equation*}
where $\sigma_{t}^2 := L^2\norms{w^{t} - w^{\star}}^2$ if $F$ is $L$-average smooth, and $\sigma^2_{t} := 4L[ F(w^{t}) - F(w^{\star}) + F(\hat{w}^s) - F(w^{\star})]$ if $F$ is convex and $L$-average smooth.

The second estimator is SARAH \cite{nguyen2017sarah}, which is expressed as follows:
\begin{equation}\label{eq:sarah_estimator}
v^{t} := v^{t-1} + [\nabla{\mbf{F}}(w^{t}, \mcal{S}_{t}) - \nabla{\mbf{F}}(w^{t-1}, \mcal{S}_{t})].
\end{equation}
It is called a stochastic recursive gradient estimator.
Unfortunately, this estimate is biased, i.e. $\mbb{E}_{\mcal{S}_{t}}[ v^{t} ] \neq \nabla{F}(w^{t})$.
However,  one can prove that 
\begin{equation*}
\mbb{E}_{\mcal{S}_{t}}[ v^{t} ] = \nabla{F}(w^{t}) + e_{t} \quad\text{and}\quad \mbb{E}_{\mcal{S}_{t}}[ \norms{ v^{t} - \nabla{F}(w^{t})}^2 ] \leq \sigma_{t}^2, 
\end{equation*}
where $e_{t} := v^{t-1} - \nabla{F}(w^{t-1})$ is an error, and $\sigma^2_{t} \leq \sigma^2_{t-1} + \frac{L^2}{b_t}\norms{w^t - w^{t-1}}^2$ if $\mbf{F}$ is $L$-average smooth, see \cite{Pham2019}.

Another interesting estimator is the hybrid variance reduced estimator proposed in \cite{Tran-Dinh2019a}, which can be written as
\begin{equation}\label{eq:hybrid_estimator}
v^{t} := (1-\beta_{t})[v^{t-1} + [\nabla{\mbf{F}}(w^{t}, \mcal{S}_{t}) - \nabla{\mbf{F}}(w^{t-1}, \mcal{S}_{t})]] + \beta_{t} u^{t},
\end{equation}
where $\beta_{t} \in [0, 1]$ and $u^{t}$ is an unbiased estimator of $\nabla{F}(w^{t})$ with variance $\hat{\sigma}_t^2$, i.e. $\mbb{E}[\norms{u^t - \nabla{F}(w^t)}^2 \mid \mcal{F}_t] \leq \hat{\sigma}_t^2$.
Again, as proven in \cite{Tran-Dinh2019a}, this is a biased estimator of $\nabla{F}(w^t)$ and if $\mbf{F}$ is $L$-average smooth, then $v^t$ satisfies $\mbb{E}[\norms{v^t - \nabla{F}(w^t)}^2 \leq \sigma_t^2$, where
\begin{equation*}
\arraycolsep=0.2em
\begin{array}{lcl}
%\mbb{E}[\norms{v^t - \nabla{F}(w^t)}^2 \mid \mcal{F}_t] & \leq &  (1-\beta_t)\norms{v^{t-1} - \nabla{F}(w^{t-1})}^2  + 2\beta_t^2\sigma^2 \vspace{1ex}\\
%&& + {~} \frac{2(1-\beta_t)^2L^2}{b_t}\norms{w^t - w^{t-1}}^2.
\sigma_t^2 \leq (1 - \beta_t)^2\sigma_{t-1}^2 +  {~} \frac{2(1-\beta_t)^2L^2}{b_t}\norms{w^t - w^{t-1}}^2 + 2\beta_t^2\hat{\sigma}_t^2. 
\end{array}
\end{equation*}
One simple choice of $u^t$ is $u^t := \nabla{\mbf{F}}(w^{t}, \mcal{S}_{t})$.
In this case, we have $\hat{\sigma}_t = \frac{\sigma^2}{b_t}$.

\subsection{Unified convergence analysis}
Similar to Subsection~\ref{subsec:conv_analysis}, let us first present our general and unified convergence analysis approach and then illustrate it through three different methods. 

\vspace{1ex}
\noindent\textbf{(a) General approach.}
Let us identify what the crucial steps in convergence analysis of Algorithm~\ref{alg:A1} are.
One of the most important steps is to establish a recursive estimate w.r.t. inner iterations $t$ of the form \eqref{eq:recursive_est}, but in conditional expectation, i.e.:
\begin{equation}\label{eq:recursive_est2}
\mbb{E}[ D_{t+1} \mid \mcal{F}_{t} ] + \Delta_{t} \leq \omega_{t} \cdot D_{t} + E_{t},
\end{equation}
where the related quantities are defined similarly to \eqref{eq:recursive_est}.
If we take the total expectation on both sides of \eqref{eq:recursive_est2}, and assume that $\omega_t = \frac{\xi_t}{\xi_{t+1}}$ for $\xi_t > 0$ and $\mbb{E}[E_t] \leq \theta_t^2M^2$ for some $M \geq 0$ and $\theta_t > 0$, then we have 
\begin{equation*} 
\xi_{t+1}\mbb{E}[ D_{t+1} ] + \xi_{t+1}\mbb{E}[\Delta_{t} ] \leq \xi_{t} \cdot \mbb{E}[ D_{t} ] + \xi_{t+1}\theta_t^2M^2.
\end{equation*}
By induction, we have 
\begin{equation}\label{eq:sgd_rates} 
\xi_{T+1}\mbb{E}[ D_{T+1} ] +  \sum_{t=0}^T\xi_{t+1}\mbb{E}[\Delta_t] \leq \xi_0\mbb{E}[D_0]  +  M^2\sum_{t=0}^T\xi_{t+1}\theta_t^2.
\end{equation}
Let $S_T := \sum_{t=0}^T\gamma_t$ with given weights $\gamma_t > 0$ (usually depending on $\xi_t$ and/or $\theta_t$).
Dividing both sides of \eqref{eq:sgd_rates} by $S_T$, we obtain
\begin{equation}\label{eq:sgd_rates2} 
\frac{1}{S_T}\sum_{t=0}^T\xi_{t+1}\mbb{E}[\Delta_t] \leq \frac{\xi_0\mbb{E}[D_0]}{S_T} + \frac{M^2}{S_T}\sum_{t=0}^T\xi_{t+1}\theta_t^2.
\end{equation}
Both estimates \eqref{eq:sgd_rates} and \eqref{eq:sgd_rates2} will allow us to estimate convergence rates of the underlying algorithm.
Let us apply this approach to prove convergence of some variants of Algorithm~\ref{alg:A1}.

\vspace{1ex}
\noindent\textbf{(b)~SGD for nonsmooth convex problems.} Let us analyze the convergence of the SGD scheme \eqref{eq:sgd_update}.
Using the update \eqref{eq:sgd_update}, we have $\norms{w^{t+1} - w^{\star}}^2  =  \norms{w^t - w^{\star}}^2 - 2\eta_t\iprods{v^t, w^t - w^{\star}} + \eta_t^2\norms{v^t}^2$.
Taking conditional expectation $\mbb{E}[\cdot \mid\mcal{F}_t]$ of this estimate and noting that $\mbb{E}[v^t \mid\mcal{F}_t] = \nabla{F}(w^t)$, we have 
\begin{equation*}
\arraycolsep=0.2em
\begin{array}{lcl}
\eta_t\iprods{\nabla{F}(w^t), w^t - w^{\star}} & = & \frac{1}{2}\norms{w^t - w^{\star}}^2 -  \frac{1}{2}\mbb{E}[\norms{w^{t+1} - w^{\star}}^2\mid\mcal{F}_t] \vspace{1ex}\\
&& + {~} \frac{\eta_t^2}{2}\mbb{E}[\norms{v^t}^2 \mid\mcal{F}_t].
\end{array}
\end{equation*}
If $F$ is convex, then we have $F(w^t) - F(w^{\star}) \leq \iprods{\nabla{F}(w^t), w^t - w^{\star}}$.
Moreover, we also have $\mbb{E}[\norms{v^t}^2 \mid\mcal{F}_t] \leq M^2$.
Combining these two expressions and the last inequality, we have
\begin{equation*}
\arraycolsep=0.2em
\begin{array}{lcl}
\underbrace{ \tfrac{1}{2}\mbb{E}[\norms{w^{t+1} - w^{\star}}^2 \mid \mcal{F}_t] }_{\mbb{E}[D_{t+1} \mid \mcal{F}_t]} + \underbrace{\eta_t[F(w^t) - F(w^{\star}) ]}_{\Delta_t} & \leq & \underbrace{\tfrac{1}{2}\norms{w^t - w^{\star}}^2}_{D_t}  + \underbrace{ \tfrac{\eta_t^2}{2}M^2}_{E_t}. %\mbb{E}[\norms{v^t}^2 \mid\mcal{F}_t].
\end{array}
\end{equation*}
This is exactly the recursive estimate \eqref{eq:recursive_est2}.
Using \eqref{eq:sgd_rates2}, we can show that
\begin{equation}\label{eq:sgd_bound3}
\arraycolsep=0.2em
\begin{array}{lcl}
\mbb{E}[F(\hat{w}) - F(w^{\star}) ] & \leq & \frac{1}{S_T}\sum_{t=0}^T \eta_t\mbb{E}[F(w^t) - F(w^{\star}) ] \vspace{1ex}\\
& \leq & \frac{1}{2S_T} \norms{w^0 - w^{\star}}^2  + \frac{M^2}{2S_T}\sum_{t=0}^T\eta_t^2,
\end{array}
\end{equation}
where $S_T := \sum_{t=0}^T\eta_t$ and $\hat{w} := \frac{1}{S_T}\sum_{t=0}^T\eta_tw^t$.
If we choose $\eta_t := \frac{C}{\sqrt{T+1}}$ for some $C > 0$, then $S_T = C\sqrt{T+1}$ and $\sum_{t=0}^T\eta_t^2 = C^2$.
In this case, \eqref{eq:sgd_bound3} becomes 
\begin{equation*}
\mbb{E}[F(\hat{w}) - F(w^{\star}) ] \leq \frac{\norms{w^0 - w^{\star}}^2}{2C\sqrt{T+1}} + \frac{M^2C}{2\sqrt{T+1}}.
\end{equation*}
If we choose $\eta_t := \frac{C}{\sqrt{t+1}}$ for some $C > 0$, then $S_T := C\sum_{t=0}^T\frac{1}{\sqrt{t+1}} \geq  2C\int_1^{T+1}\frac{1}{2\sqrt{t}}dt  = 2C(\sqrt{T+1} - 1)$ and $\sum_{t=0}^T\eta_t^2 = C^2\sum_{t=0}^T\frac{1}{t+1} \leq C^2(1 + \ln(T+1))$.
In this case,  \eqref{eq:sgd_bound3} becomes 
\begin{equation*}
\mbb{E}[F(\hat{w}) - F(w^{\star}) ] \leq \frac{\norms{w^0 - w^{\star}}^2}{4C(\sqrt{T+1} - 1)} + \frac{M^2C(1 + \ln(T+1))}{4(\sqrt{T+1} - 1)}.
\end{equation*}

\vspace{1ex}
\noindent\textbf{(c)~SGD for smooth and nonconvex problems.}
We consider the case $F$ is $L$-smooth.
In addition, we assume that our stochastic estimator $v^t$ is unbiased, i.e. $\mbb{E}[v^t \mid \mcal{F}_t] = \nabla{F}(w^t)$ and has bounded variance as $\mbb{E}[\norms{v^t - \nabla{F}(w^t)}^2 \mid \mcal{F}_t] \leq \sigma^2$.
In this case, we have $\mbb{E}[ \norms{v^t}^2 \mid\mcal{F}_t] \leq \norms{\nabla{F}(w^t)}^2 + \sigma^2$.
Using this inequality,  $\mbb{E}[v^t \mid \mcal{F}_t] = \nabla{F}(w^t)$, 
%\textcolor{red}{Removed $\mbb{E}[v^t - \nabla{F}(w^t) \mid\mcal{F}_t] = 0$: From $\mbb{E}[v^t \mid \mcal{F}_t] = \nabla{F}(w^t)$ we can already infer $\eta_t \norms{\nabla{F}(w^t)}^2 = \eta_t\mbb{E}[ \iprods{\nabla{F}(w^t), v^t  } \mid \mcal{F}_t]$. Hence, we can shorten the formula after the second inequality below and save a line.} 
and  the $L$-smoothness of $F$, we can derive 
\begin{equation*} 
\arraycolsep=0.2em
\begin{array}{lcl}
\mbb{E}[F(w^{t+1}) \mid \mcal{F}_t] &\leq& F(w^t) - \eta_t\mbb{E}[\iprods{\nabla{F}(w^t), v^t} \mid \mcal{F}_t] + \frac{L\eta_t^2}{2}\mbb{E}[ \norms{v^t}^2 \mid \mcal{F}_t] \vspace{1ex}\\
&\leq & F(w^t) - \eta_t\norms{\nabla{F}(w^t)}^2 
%- \eta_t\mbb{E}[ \iprods{\nabla{F}(w^t), v^t - \nabla{F}(w^t)} \mid \mcal{F}_t] \vspace{1ex}\\
%&& 
+ {~} \frac{L\eta_t^2}{2}\norms{\nabla{F}(w^t)}^2 + \frac{L\eta_t^2\sigma^2}{2} \vspace{1ex}\\
&= & F(w^t) - \eta_t\big(1 - \tfrac{L\eta_t}{2}\big)\norms{\nabla{F}(w^t)}^2 +  \frac{L\eta_t^2\sigma^2}{2}.
\end{array}
\end{equation*}
This inequality leads to
\begin{equation*} 
\arraycolsep=0.2em
\begin{array}{lcl}
\mbb{E}[\underbrace{F(w^{t+1}) - F^{\star}}_{D_{t+1}} \mid \mcal{F}_t] + \underbrace{\eta_t\big(1 - \tfrac{L\eta_t}{2}\big)\norms{\nabla{F}(w^t)}^2}_{\Delta_t} \leq \underbrace{F(w^t) - F^{\star}}_{D_t} + \underbrace{\tfrac{L\eta_t^2\sigma^2}{2}}_{E_t},
\end{array}
\end{equation*}
which is exactly \eqref{eq:recursive_est2} with $\omega_t = 1$, provided that $0 < \eta_t < \frac{2}{L}$.
By using this estimate we can derive a convergence rate for $\frac{1}{S_T}\sum_{t=0}^T\gamma_t\mbb{E}[\norms{\nabla{F}(w^t)}^2]$ with $\gamma_t := \eta_t\big(1 - \tfrac{L\eta_t}{2}\big)$ and $S_T := \sum_{t=0}^T\gamma_t$ as done in \cite{ghadimi2013stochastic}.
We omit the details here.
%Since $\mbb{E}[\norms{v^t}^2 \mid\mcal{F}_t] \leq M^2$% = \mbb{E}[\norms{v^t - \nabla{F}(w^t)}^2 \mid\mcal{F}_t] + \norms{\nabla{F}(w^t)}^2 \leq \sigma^2 + \norms{\nabla{F}(w^t)}^2$.

\vspace{1ex}
\noindent\textbf{(d)~Hybrid variance-reduced SGD for smooth and nonconvex problems.}
We analyze one variance-reduced variant of Algorithm~\ref{alg:A1} where the inner loop updates $w^{t+1} := w^t - \eta_tv^t$ with $v^t$ being given by \eqref{eq:hybrid_estimator} for $b_t = 1$, see \cite{Tran-Dinh2019a}.
In addition, we do not need the outer loop, leading to a \textbf{single-loop} algorithm.

Let us analyze its convergence rate.
First, by the $L$-smoothness of $F$ and the relation 
%\textcolor{red}{This allowed me to remove the second inequality in the derivation below.} 
$-2\iprods{a, b} = \norms{a-b}^2 - \norms{a}^2 - \norms{b}^2$, we can derive
\begin{equation*} 
\arraycolsep=0.2em
\begin{array}{lcl}
\mbb{E}[F(w^{t+1}) \mid \mcal{F}_t] &\leq& F(w^t) - \eta_t\mbb{E}[\iprods{\nabla{F}(w^t), v^t} \mid \mcal{F}_t] + \frac{L\eta_t^2}{2}\mbb{E}[ \norms{v^t}^2 \mid \mcal{F}_t] \vspace{1ex}\\
%&\leq & F(w^t) - \eta_t\norms{\nabla{F}(w^t)}^2 - \eta_t\mbb{E}[ \iprods{\nabla{F}(w^t), v^t - \nabla{F}(w^t)} \mid \mcal{F}_t] \vspace{1ex}\\
%&& + {~}  \frac{L\eta_t^2}{2}\mbb{E}[ \norms{v^t}^2 \mid \mcal{F}_t] \vspace{1ex}\\
&= & F(w^t) - \tfrac{\eta_t}{2}\norms{\nabla{F}(w^t)}^2 +  \frac{\eta_t}{2}\mbb{E}[ \norms{v^t - \nabla{F}(w^t)}^2 \mid \mcal{F}_t] \vspace{1ex}\\
&& - {~}  \tfrac{\eta_t}{2}(1 - L\eta_t)\mbb{E}[ \norms{v^t}^2 \mid \mcal{F}_t].
\end{array}
\end{equation*}
Since $\mbb{E}[ \norms{v^t - \nabla{F}(w^t)}^2 \mid \mcal{F}_t]  \leq \sigma_t^2$ and $0 < \eta_t \leq \frac{1}{L}$, this inequality reduces to
\begin{equation*} 
\arraycolsep=0.2em
\begin{array}{lcl}
\mbb{E}[F(w^{t+1}) - F^{\star} + \tfrac{\eta_t(1 - L\eta_t)}{2} \norms{v^t}^2  \mid \mcal{F}_t] &\leq  F(w^t)  - F^{\star} - \tfrac{\eta_t}{2}\norms{\nabla{F}(w^t)}^2 +  \frac{\eta_t\sigma_t^2}{2}.
\end{array}
\end{equation*}
Since $\sigma_t^2 \leq (1-\beta_t)^2\sigma_{t-1}^2 + 2(1-\beta_t)^2L^2\norms{w^t - w^{t-1}}^2 + 2\beta_t^2\hat{\sigma}^2$ and $w^t - w^{t-1} = -\eta_{t-1}v^{t-1}$, we have
\begin{equation*}
\sigma_t^2 \leq (1-\beta_t)^2\sigma_{t-1}^2 + 2L^2(1-\beta_t)^2\eta_{t-1}^2\norms{v^{t-1}}^2 + 2\beta_t^2\hat{\sigma}^2.
\end{equation*}
Multiplying this inequality by $\frac{c_t}{2} > 0$ and adding to the last estimate, we obtain
\begin{equation*} 
\arraycolsep=0.1em
\begin{array}{lcl}
\mbb{E}[F(w^{t+1}) - F^{\star} & + &  \tfrac{\eta_t(1 - L\eta_t)}{2} \norms{v^t}^2 \mid \mcal{F}_t] + \frac{(c_t-\eta_t)}{2}\sigma_t^2 + \tfrac{\eta_t}{2}\norms{\nabla{F}(w^t)}^2 \leq F(w^t) - F^{\star} \vspace{1ex}\\
&& + {~}  L^2c_t(1-\beta_t)^2\eta_{t-1}^2\norms{v^{t-1}}^2 + \frac{c_t(1-\beta_t)^2}{2}\sigma_{t-1}^2  + c_t\beta_t^2\hat{\sigma}^2.
\end{array}
\end{equation*}
For simplicity, we choose all parameters to be constant.
Let us define $D_t := F(w^t) - F^{\star} + \frac{\eta(1-L\eta)}{2}\norms{v^{t-1}}^2 + \frac{(c -\eta)}{2}\sigma_{t-1}^2$, and impose the following conditions:
\begin{equation}\label{eq:para_conds}
2L^2\eta^2c(1-\beta)^2 \leq \eta(1-L\eta)\quad\text{and} \quad c(1-\beta)^2 \leq c - \eta.
\end{equation}
Then, the last estimate leads to 
\begin{equation*} 
\mbb{E}[ D_{t+1} \mid \mcal{F}_t]  + \underbrace{\tfrac{\eta}{2}\norms{\nabla{F}(w^t)}^2}_{\Delta_t} \leq D_t  + \underbrace{c\beta^2\hat{\sigma}^2}_{E_t},
\vspace{-1ex}
\end{equation*}
which is exactly \eqref{eq:recursive_est2} with $\omega_t = 1$.

Assume that we choose $\eta \in (0, \frac{1}{L})$ and $c > 0$ such that $2L^2\eta^2(c-\eta) = \eta(1-L\eta)$, leading to $c := \frac{1-L\eta}{2L^2\eta} + \eta = \frac{1-L\eta + 2L^2\eta^2}{2L^2\eta}$.
Moreover, $(1- \beta)^2 \leq 1-  \frac{2L^2\eta^2}{1 - L\eta + 2L^2\eta^2}$.
Then, both conditions of \eqref{eq:para_conds} hold with equality. %s for $(1-\beta)^2\leq \frac{1-L\eta}{1-L\eta + 2L^2\eta^2}$.
In this case, we obtain $\mbb{E}[ D_{t+1} \mid \mcal{F}_t]  + \tfrac{\eta}{2}\norms{\nabla{F}(w^t)}^2 \leq D_t  + \tfrac{(1-L\eta + 2L^2\eta^2)\beta^2}{2L^2\eta}\hat{\sigma}^2$.
This inequality implies 
%\textcolor{red}{In the second inequality we need $\sigma_{-1}$ rather than $\sigma_0$ -- this goes back to the recursive inequality below (1.25), what can we choose for $\sigma_{-1}$?}
%\textcolor{red}{We also need $1-L\eta + 2L^2\eta^2\leq 1$, that is, $\eta\leq \frac{1}{2L}$, a factor 2 smaller compared to what was assumed before.}
\begin{equation*} 
\arraycolsep=0.1em
\begin{array}{lcl}
\tfrac{1}{T+1}\sum_{t=0}^T\mbb{E}[\norms{\nabla{F}(w^t)}^2] & \leq & \tfrac{2}{\eta(T+1)}D_0 + \tfrac{(1-L\eta + 2L^2\eta^2)\beta^2}{L^2\eta^2}\hat{\sigma}^2 \vspace{1ex}\\
&\leq & \tfrac{2[F(w^0) - F^{\star}]}{\eta(T+1)} +   \tfrac{\norms{v^0}^2}{(T+1)} +  \tfrac{\sigma_{-1}^2}{2L^2\eta^2(T+1)} + \tfrac{\beta^2\hat{\sigma}^2}{L^2\eta^2}.
\end{array}
\end{equation*}
%We also need $1-L\eta + 2L^2\eta^2\leq 1$, that is, $\eta\leq \frac{1}{2L}$, a factor 2 smaller compared to what was assumed before.
Finally, we choose $\eta := \frac{1}{L(T+1)^{1/3}} \leq \frac{1}{L}$, $\sigma_{-1} := \frac{1}{(T+1)^{1/3}}$, and 
%\textcolor{red}{Does this fit the extra condition on $\beta$ in red above (may be $\beta$ needs another factor $\sqrt{2}\cdot L$)?} 
$\beta := \BigO{\frac{1}{(T+1)^{2/3}}}$ such that $(1- \beta)^2 \leq 1-  \frac{2L^2\eta^2}{1 - L\eta + 2L^2\eta^2}$ (always exist such a $\beta$).
Moreover, the last estimate shows that 
\begin{equation*}
\tfrac{1}{T+1}\sum_{t=0}^T\mbb{E}[\norms{\nabla{F}(w^t)}^2] = \BigO{\frac{1}{(T+1)^{2/3}}}, 
\end{equation*}
as proven in  \cite{Tran-Dinh2019a}.

We have illustrated our approach by using the recursive estimate \eqref{eq:recursive_est2} to analyze the convergence of three SGD schemes, including variance-reduced methods.
We believe that this approach can be used to analyze other variants including SVRG and SARAH.

\section{Concluding remarks}\label{chapter_XX_sec_conclusion}
We have reviewed several main components that constitute the gradient descent method and its variants, including deterministic and stochastic ones, ranging from convex to nonconvex problems.
We have provided a simple and unified convergence analysis framework relying on an elementary recursive estimate under the most basic structure assumptions commonly used in the literature. 
While this approach can be applied to analyze several methods, we have only illustrated it on a few well-known schemes.
Note that we have not proposed any new algorithms, but rather unified the convergence analysis using a simple recursive estimate.
However, we believe that such an approach can be extended beyond what we have done in this paper.
The following research topics are interesting to us.
First, can one still apply our analysis to accelerated variance-reduced stochastic gradient-type methods?
Perhaps, this can possibly be done by using the idea from a recent work \cite{driggs2019accelerating}.
Second, how can we extend our framework to study other optimization methods in distributed systems and federated learning?
We emphasize that many algorithms in these fields can be viewed as a randomized [block-]coordinate methods.
Therefore, extensions to coordinate methods and shuffling methods are promising and remain open.
Third, is it possible to extend and adapt our analysis to asynchronous gradient-based algorithms?
We believe that such an extension is possible as long as the delay is bounded.
However, one needs to modify the recursive expression to capture with the delayed updates, leading to an extra error term in the recursive inequality. 
Finally, our approach can be used to analyze convergence of algorithms for minimax and variational inequality problems, which have recently gained tremendous popularity \cite{Facchinei2003,TranDinh2020f}.

%\textcolor{red}{Do you like to add the problem of extending the unified framework to include asynchronous updates on the weight vector? There may not be enough space to make the reader understand what we mean by this: We would want to explain that at each entry in the weight vector we add a computed  "update" within a "delay" from the moment the update  was computed. This reflects entry by entry asynchronous behavior where write and read from entries are assumed atomic.}

\section{Acknowledgements}\label{chapter_XX_sec_acknowledgement}
The work of Q. Tran-Dinh is partly supported by the Office of Naval Research [grant number ONR-N00014-20-1-2088] (2020–2023) and the National Science Foundation (NSF) [grant number NSF DMS-2134107] (2022-2027). 

%\section{References}

%\bibliography{reference}
\bibliographystyle{plain}
%\bibliography{/Users/quoctd/Dropbox/E-Books/tran_bibtex_new}

% \begin{figure}
% \title{}% Figure title
% \caption{}% Figure caption
% \label{}%
% \end{figure}

% \begin{table}
% \title{}% Table title
% \caption{}% Table caption
% \label{}%
% \end{table}
\end{document}